\documentclass[12pt]{article}\usepackage{amssymb,amsmath}
\usepackage{amscd}

\usepackage{times}
\usepackage{amssymb,amsmath}
\usepackage{latexsym}
\usepackage{amsthm}
\usepackage{enumerate}

\newtheorem{theorem}{Theorem}[section]
\newtheorem{lemma}[theorem]{Lemma}
\newtheorem{corollary}[theorem]{Corollary}
\newtheorem{proposition}[theorem]{Proposition}

\newfont{\Goth}{eufm10 at 11pt}

\newfont{\biggoth}{eufm10 at 18pt}

\newcommand{\HH}{{\mathcal H}}
\newcommand{\Q}{\mathbb Q}
\newcommand{\mgoth}[1]{\mbox{\Goth #1}}

\newcommand{\R}{\mathbb R}
\newcommand{\Z}{\mathbb Z}

 \newcommand{\heis}{\mathcal Heis}
\newcommand{\pr}{\prime}

\author{ Andr\'e DIATTA \and Brendan FOREMAN
\footnote{ \footnotesize University of Liverpool. Department of Mathematical Sciences. M$\&$O Building, Peach Street, Liverpool, L69 7ZL, UK. adiatta@liv.ac.uk.
\newline
\noindent
$\bullet$ A part of this work was done while the author was supported by the IST Programme of the European Union (IST-2001-35443).}}

\title{\bf  Lattices in contact Lie groups and 5-dimensional contact solvmaniolds}

\begin{document}

\maketitle

\begin{abstract}{ \footnotesize
We investigate the existence and  properties of uniform lattices in Lie
 groups and use these results to prove that, in dimension 5, there are
 exactly seven connected and simply connected contact Lie groups with
 uniform lattices, all of which are solvable.  Issues of symplectic
 boundaries are explored, as well. It is also shown that the special affine
 group has no uniform lattice.
}
\end{abstract}


\section{Introduction}\label{chap:introduction}

This paper investigates the geometry of compact contact manifolds that are uniformized by contact Lie groups, i.e., manifolds of the form $\Gamma\setminus G$ for some Lie group $G$ with a left invariant contact structure and uniform lattice $\Gamma\subset G$. In particular, we restrict our attention to dimension five and describe which five-dimensional contact Lie groups admit uniform lattices.  We prove that there are exactly seven connected and simply connected such Lie groups.  Five of them are central extensions; the other two are semi-direct products.  Furthermore, all seven are solvable. Let us remind that, in the symplectic counterpart, there are only $4$ connected and simply connected Lie groups with a  lattice, that can bear a left invariant symplectic form \cite{medina-revoy-lattice07}.

This paper is organized as follows. In Section
\ref{chap:preliminaries}, we give the preliminaries for the work ahead.  This includes both a review of several classical results and some original results regarding contact Lie groups.  Fundamental to this paper are Theorem \ref{thm:Diatta1}, which describes all five-dimensional contact Lie algebras, and the
list in Subsection \ref{5-dimensional-unimod}, which delineates the Lie
algebras of all the five-dimensional unimodular contact Lie groups.

In
Section \ref{chap:main results}, the main theorem of the paper (Theorem
\ref{thm:listing}) is stated as well as an immediate corollary.  This
theorem is proven in Section \ref{chap:big proof}.  A major yet technical aspect of this proof is the list of certain structures on the Lie algebras of the Lie groups in Subsection \ref{5-dimensional-unimod}.  For ease of reading, this list has been relegated to Appendix I (Section \ref{nillist}).

Finally, Section \ref{chap:symplectic-disconnected-boundary} constructs compact symplectic ($2n+2$)-manifolds whose boundaries are disconnected contact ($2n+1$)-manifolds uniformized by contact Lie groups and hence, when $n=2,$ by the Lie groups of Theorem \ref{thm:listing}.  This is a generalisation to all higher dimensions of a construction used in
 \cite{Geiges95}, to give counter-examples, when $n=1$, to the question of E. Calabi as to whether symplectic compact manifolds with a boundary of contact type, admit a connected boundary, as it is the case for compact complex manifolds with  strictly pseudo-convex boundary. The counterexamples in \cite{Geiges95} encompass those by D. McDuff in  \cite{mcduff}.

The authors thank the referee of this paper for the many essential hints and suggestions.  This paper would not have been possible without this help.


\section{Preliminaries} \label{chap:preliminaries}

\subsection{Lattices on solvable Lie groups}\label{subsec:solvlattices}

A {\it lattice} of a Lie group $G$ is a discrete subgroup
$\Gamma$ such that the manifold $\Gamma\setminus G$ has a finite volume.
If $\Gamma\setminus G$ is compact, then $\Gamma$ is called a {\it uniform lattice}.
This is a well-studied field, and much of the following material has been derived from Chapter 2 of Part I in \cite{ov2}, much of which is itself an exposition of classical results by Mostow (\cite{mostow}), Auslander (\cite{aus1},\ \cite{aus2}) and Raghunathan (\cite{raghu}).  More details as well as more results on this topic can be found within these various sources.

One of the most important results on lattices on general Lie groups was proved by Milnor in \cite{milnor}.
\begin{theorem}
If $G$ is a Lie group with a uniform lattice, then its Lie algebra is unimodular.
\end{theorem}
\noindent For nilpotent Lie groups, more precise results are known.  In particular, a lattice on a nilpotent Lie group induces lattices on the central series of the Lie group.
\begin{theorem}\label{handy}
Let $N$ be a simply connected nilpotent Lie group with lattice $\Gamma$.  Let
$\dots \subset N_2 \subset N_1 \subset N_0=N$ be the decreasing central series of $N$.  
Then, for each $j=0,\ 1,\ 2,\ \dots,$ $\Gamma\cap N_j$ is a lattice of $N_j.$
\end{theorem}
\noindent And, most importantly, there is a well-known necessary and sufficient condition for the existance of a lattice on a given nilpotent Lie group.
\begin{theorem}
Let $N$ be a nilpotent Lie group.  Then $N$ has a lattice if and only if
its Lie algebra $\mgoth{n}$ has a $\Q$-algebra $\mgoth{n}_\Q$, that is, $\mgoth{n}$ 
has a basis whose Lie structure constants are integers.
\end{theorem}

For solvable Lie groups, several general results are well known.  In order to describe these results, we first review some structure theory on solvable Lie groups.
Let $G$ be a simply-connected solvable Lie group with Lie algebra $\mgoth{g}.$  Let $N$ be the nilradical of $G$ with corresponding Lie algebra $\mgoth{n}$, i.e., $N$ is the 
maximal nilpotent normal subgroup of $G$ so that $\mgoth{n}$ is the maximal nilpotent
ideal of $\mgoth{g}$.  This induces a short exact sequence
$$1 \rightarrow N \rightarrow G \rightarrow T \rightarrow 1,$$
where $T$ is the Abelian group given by $T=N\setminus G.$
$G$ is called {\it splittable}, if this sequence splits, that is, there is a right inverse homomorphism of the projection $G \rightarrow T.$  This condition is equivalent to 
the existence of a homomorphism $b: T\rightarrow Aut(N)$ such that $G$ is 
isomorphic to the semi-direct product $N \rtimes_b T.$

The first known result regarding lattices of solvable groups indicates just how much more specialized subgroups lattices are for solvable Lie groups than they are for general Lie groups.
\begin{theorem}
A lattice on a solvable Lie group is a uniform lattice.
\end{theorem}
\noindent Furthermore, we have Mostow's well-known result.
\begin{theorem}[Mostow \cite{mostow}]  \label{thm:mostow} Let $\Gamma$ be a lattice in a connected solvable Lie group $G$ with nilradical $N.$  Then $\Gamma\cap N$ is a lattice of $N.$
\end{theorem}

In \cite{wang}, Wang proved the general structure of a lattice of a solvable Lie group.
\begin{theorem}[Wang \cite{wang}]\label{thm:wang}
A group $\Gamma$ is isomorphic to a discrete subgroup in a simply-connected Lie group
if and only if there is a lattice $\Delta$ of a simply-connected nilpotent Lie group 
and non-negative integer $k$ such that
$$0 \rightarrow \Delta \rightarrow \Gamma \rightarrow \Z^k \rightarrow 0$$
is a short, exact sequence.
\end{theorem}
In particular, this implies that, if $G=N \rtimes_b T$ is a simply-connected splittable solvable Lie group
with nilradical $N$ and $\Gamma$ is a lattice of $G$, then 
$\Gamma$ is isomorphic to $\Delta \rtimes_b T_\Z$ where 
$\Delta$ is a lattice of $N$ and $T_\Z$ a lattice of $T$ such that $b(T_\Z)\subset Aut(\Delta).$

\subsection{Heisenberg groups}\label{Heisenberg}

 Besides $\R^m$ under addition, the most encountered Lie group in the work below will be the Heisenberg groups in three and five dimensions.  In general, the {\it $(2n+1)$-dimensional Heisenberg group} $\heis^{2n+1}$ is the subgroup of $Sl(n+2,\R)$ given by 
$$ \heis^{2n+1}=
\left\{\sigma=\begin{pmatrix}1&y&x\\0&I_{n}&z^t\\0&0&1\\ \end{pmatrix}:
y, z\in \R^n, ~x\in\R\right\},$$ 
where the column vector $z^t$ is the transpose
of the vector $z=(z_1,\ldots, z_n)$ and $I_{n}$ is the identity map of
$\R^{n}.$  Equivalently, $\heis^{2n+1}$ can be considered as the central extension of the symplectic Lie group
 $\R^{2n}$ under addition with the standard symplectic form $\omega_1$.

 The Lie algebra of $\heis^{2n+1}$ is given by
 $$\mgoth{h}_{2n+1} =\left\{ X=\begin{pmatrix}0&b&a\\0&0&c^t\\0&0&0\\
				   \end{pmatrix}: b, c\in \R^n,
 a\in\R\right\}.$$ 
For $i,\ j\in\{1,\dots, n+2\},$ let $e_{i,j}$ be the $(n+2)\times (n+2)$
 matrix, all of whose entries are zero except the ${ij}$-th entry which is
 equal to $1.$ We set $e_1:=e_{1,n+2}$, $e_k:=e_{1,k}$ and
 $e_{n+k}:=e_{k,n+2}$ for $k=2, \dots, n+1.$  Then $\{e_1, \dots, e_{2n+1} \}$ is a basis of $\mgoth{h}_{2n+1}$ with exactly $n$ nontrivial Lie brackets
 relations, namely, $[e_{k},e_{n+k}]=e_1$ for all $k=2,\dots, n+1.$ 
 If we let
 $(e_1^*,\dots, e_{2n+1}^*)$ stand for the dual basis of
 $(e_1,\dots, e_{2n+1})$, then $e_1^*$ is a contact form on  $\mgoth{h}_{2n+1}.$  In terms of the original coordinates on $\heis^{2n+1},$ the left invariant vector fields are given by 
$e_1^+= \frac{\partial}{\partial x}$, $e_k^+= \frac{\partial}{\partial y_{k-1}},$ $e_{n+k}^+= \frac{\partial}{\partial z_{k-1}}+y_{k-1}\frac{\partial}{\partial y}$ for $k=2,\dots, n+1.$ The left invariant contact form on $\heis^{2n+1}$
 corresponding to $e_1^*$ is $e_1^{*,+}=dx-\displaystyle \sum_{i=1}^ny_idz_i.$

 The exponential map $\exp:\mgoth{h}_{2n+1}\rightarrow \heis^{2n+1}$
 is a diffeomorphism, and we denote its inverse by $\log.$ Specifically, these mappings are given by
 
 \begin{eqnarray*}
 \exp \begin{pmatrix}0&a&c\\0&0_n&b^t\\0&0&0\\
				   \end{pmatrix} &=& \begin{pmatrix}1&a&c+\frac{1}{2}ab^t\\ 0&I_n&b^t\\ 0&0&1 \end{pmatrix},\\  
 \log \begin{pmatrix}1&x&z\\0&I_{n}&y^t\\0&0&1\\ \end{pmatrix} &=& \begin{pmatrix} 0&x&z-\frac{1}{2}xy^t\\0&0_n&y^t\\0&0&0 \end{pmatrix}.\\ \end{eqnarray*}

We focus in on the case where $n=1.$
Let $N$ be the Lie group given by $N=\heis^3 \times \R$.  Its Lie algebra is given by
$\mgoth{n}=\mgoth{h}_3\oplus \R$.  
Furthermore $N$ has a left
 invariant nondegenerate closed 2-form, hence defining a left invariant
 symplectic structure  $\omega = dx\wedge dz + dw\wedge dy,$ where $w$
 is the coordinate in $\R.$ It corresponds to the symplectic form
 $\omega_2={e_1}^*\wedge{e_3}^*+{e_4}^*\wedge{e_2}^*,$ on  $\mgoth{h}_{3}
 \oplus \R$,  where $\mgoth{h}_3=<e_1,e_2,e_3>_{\R}$ as above.

Since lattices of $\heis^3$ exist (see, for example, \cite{gw}), lattices on $N$ exist.  In fact,
if $\Gamma$ is a lattice of $N$, then
$$[\Gamma, \Gamma]\subset \left\{ 
\begin{pmatrix} x\\ 0\\ 0\\ 0\\ \end{pmatrix}: x\in\R\right\} \subset N.$$
In particular, $[\Gamma, \Gamma]$ is a subgroup of $\Gamma.$ So, there is $x_0\in \R^+$ such that
$$[\Gamma, \Gamma]\subset \left\{ 
\begin{pmatrix} kx_0\\ 0\\ 0\\ 0\\ \end{pmatrix}: k\in\Z\right\} \subset N.$$
We will make extensive use of this fact when we are proving that certain Lie groups do not have
lattices.

\subsection{Five-dimensional contact Lie groups}

A {\it contact Lie group} is a $(2n+1)$-dimensional Lie group $G$ with a left-invariant differential form $\eta$ such that $\eta\wedge d\eta^n\neq 0$.  Set $\HH=ker\ \eta$.  Then
$\HH$ is a left-invariant $2n$-dimensional subbundle of $TG$ so that $\HH$ induces a subspace of the Lie algebra $\mgoth{g}$ of $G$, which we will also denote as $\HH.$  An element $X\in\mgoth{g}$ is called {\it horizontal}, if $X\in\HH.$  A submanifold of $G$ is called {\it totally isotropic}, if its tangent space in $G$ is horizontal everywhere.  A totally isotropic submanifold of (maximal) dimension $n$ is called a {\it Legendrian} submanifold of $G.$

Let $\xi$ be the unique left-invariant vector field in $\mgoth{g}$ defined by
\begin{eqnarray*}
\eta(\xi)&=&1,\\
d\eta(\xi, *)&=&0.\\
\end{eqnarray*}
Then $\mgoth{g}=\left< \xi \right>_\R \oplus \HH,$ and $\xi$ is called the {\it Reeb
vector field} of $\eta.$

\begin{lemma}
Let $(G,\ \eta)$ be a solvable contact Lie group with nilradical $N$.  Let $\mgoth{n}$ be the
Lie algebra of $N$.  Then $\mgoth{n}$ is not contained in $\HH.$
\end{lemma}
\noindent {Proof:\ \ }$[\mgoth{g}, \mgoth{g}]\subset \mgoth{n}.$  Thus, if $\mgoth{n} \subset \HH,$
$[X,Y]\in \mgoth{n} \subset \HH$ for any $X,\ Y.$
If this were the case, then $d\eta(X,Y)=-\frac{1}{2} \eta([X,Y])=0$ for any $X,\ Y.$
Thus, $d\eta=0$ on $G,$ a contradiction. \qed

A Lie algebra  $\mgoth{g}$ is said to be {\it decomposable} if
it is the direct sum $\mgoth{g}= \mgoth{g}_1 \oplus \mgoth{g}_2$ of
two ideals  $\mgoth{g}_1$ and $\mgoth{g}_2.$ Such a Lie algebra has
a contact form if and only if  $\mgoth{g}_1$ has a contact form and
$\mgoth{g}_2$ an exact symplectic form, or vice versa. 

\begin{lemma}\label{lem:decomp} If a  contact Lie algebra (resp. group) is unimodular, then it is necessarily nondecomposable.\end{lemma}

\proof  A decomposable Lie algebra $\mgoth{g}= \mgoth{g}_1
\oplus \mgoth{g}_2$ is unimodular if and only if both $\mgoth{g}_1$
and $\mgoth{g}_2$ are unimodular. But as noted above, if $\mgoth{g}$
had a contact form, then  $\mgoth{g}_1$ (or $\mgoth{g}_2$)
would have an exact symplectic form. And due to the existence of a
left invariant radiant vector field for the associated left
invariant affine connection, a Lie group with a left invariant exact
symplectic form cannot be unimodular (see
\cite{diatta-medina}). This, applied to any Lie group with Lie algebra  $\mgoth{g}_1$,
would lead to a contradiction.
\qed

\begin{corollary}
Let $G=N\rtimes_b T $ be a $(2n+1)$-dimensional, simply-connected splittable solvable Lie group with nilradical $N$ and homomorphism $b: T \rightarrow Aut(N).$  Let $\eta\in\mgoth{g}^*$ be a left-invariant contact structure on $G$. Then
\begin{enumerate}
\item The subspace $\mgoth{n}\cap \HH$ has codimension 1 in $\mgoth{n}.$
\item $dim\ T\leq n$ and $dim\ \mgoth{n}\geq n+1.$
\item For every $X\in T$, there is an $X^\prime\in \mgoth{n}\cap\HH$ such that
$d\eta(X,X^\prime)=1.$
\end{enumerate}
\end{corollary}

\subsubsection{Five-dimensional solvable contact Lie algebras}

In \cite{Diatta-Contact}, the first author classified the five-dimensional simply connected contact Lie groups (via their Lie algebras) with the following theorem.

\begin{theorem}[Diatta \cite{Diatta-Contact}]\label{thm:Diatta1}Let $G$ be a five-dimensional Lie group with Lie algebra $\mgoth{g}$.

\begin{enumerate}
\item Suppose $G$ is non-solvable.  Then $G$ is a contact Lie group if and only if $\mgoth{g}$ is one of the following Lie algebras:
\begin{enumerate}
\item $aff(\R)\oplus \mgoth{sl}(2,\R)$, $aff(\R)\oplus \mgoth{so}(3,\R)$ (decomposable cases) or
\item $\mgoth{sl}(2,\R)\ltimes \R^2$ (non-decomposable case).
\end{enumerate}
\item Suppose that $G$ is solvable such that $\mgoth{g}$ is non-decomposable with trivial center $Z(\mgoth{g}).$
Then
\begin{enumerate}
\item If the derived ideal $ [\mgoth{g},\mgoth{g}]$ has dimension three
      and is non-Abelian, then $\mgoth{g}$ is a contact Lie algebra.
\item If $ [\mgoth{g},\mgoth{g}]$ has dimension four, then $\mgoth{g}$ is contact if and only if
\begin{enumerate}
\item $dim(Z([\mgoth{g},\mgoth{g}]))=1$\ or
\item $dim(Z([\mgoth{g},\mgoth{g}]))=2$ and there is a $v\in\mgoth{g}$ such that $Z([\mgoth{g},\mgoth{g}])$ is not an eigenspace of $ad_v.$
\end{enumerate}
\end{enumerate}
\end{enumerate}
\end{theorem}
\noindent  The first statement of this result taken with Lemma \ref{lem:decomp} implies that the only unimodular non-solvable five-dimensional contact Lie group is
$\mgoth{sl}(2,\R)\ltimes \R^2.$  Furthermore, the second statement
in conjunction with the list of five-dimensional solvable Lie algebras in \cite{blz} yields the list of all five-dimensional solvable contact Lie algebras, a total of 24 distinct nondecomposable Lie algebras and families of Lie algebras.  Among these, exactly 12 are unimodular.  They are listed below along with an example of a contact form $\eta$.  The label for each Lie algebra refers to that algebra's position in the original list in \cite {Diatta-Contact} and will serve as the name of that Lie algebra (or corresponding simply connectd Lie group) throughout this paper.

\subsubsection{Five-dimensional unimodular solvable contact Lie algebras} \label{5-dimensional-unimod}

\medskip
Below is the list of unimodular solvable contact Lie algebras of dimension $5$.

\medskip
\noindent{\bf Central extensions}
 \begin{description}
\item [D1]$[e_2,e_4] = e_1$, $[e_3,e_5] = e_1$, $\eta:= e_1^*$. This
is the
  Heisenberg Lie algebra $\mgoth{h}_5$. See Section \ref{Heisenberg}.
\item [D2] $[e_3,e_4] = e_1$, $[e_2,e_5] = e_1$, $[e_3,e_5] = e_2$,   $\eta:
=e^*_1$. This is the central extension
 $\mgoth{b}\times_\omega \mathbb R e_1$, where  $\omega=e_3^*\wedge
	    e_4^* + e_2^*\wedge e_5^*$ and $\mgoth{b}= \mgoth{h}_3\oplus \mathbb Re_4$, as in  \ref{Heisenberg}.
\item [D3] $[e_3,e_4] = e_1$, $[e_2,e_5] =  e_1$, $[e_3,e_5] = e_2$, $[e_4,e_5] = e_3$,
$\eta=e^*_1$. This is the central extension $\mgoth{b}\times_\omega
\mathbb R e_1$, where  $\omega=e_3^*\wedge e_4^* + e_2^*\wedge
e_5^*$ and $\mgoth{b}= span (e_2,e_3,e_4,e_5 )$ with Lie bracket
$[e_3,e_5] = e_2$, $[e_4,e_5] =e_3.$
\item [D5] \ $[e_2,e_3] = e_1$, $[e_2,e_5] =e_2$, $[e_3,e_5] = - e_3$, $ [e_4,e_5] =  e_1,$ $\eta=e^*_1$.
 This is $\mgoth{b}\times_\omega \mathbb R e_1$, where  $\omega=e_2^*\wedge e_3^* + e_4^*\wedge e_5^*$ and $\mgoth{b}= span (e_2,e_3,e_4,e_5 )$ with Lie bracket $[e_2,e_5] =e_2$, $[e_3,e_5] = - e_3.$
\item [D11] \ $ [e_2,e_3]= e_1; [e_2,e_5]= e_3; [e_3,e_5]= -e_2; [e_4,e_5]= \epsilon e_1$;
 $\epsilon =\pm 1$; $\eta=e^*_1$. Here $\mgoth{g}=\mgoth{b}\times_\omega \mathbb R e_1$,
 where  $\omega=e_2^*\wedge e_3^* + e_4^*\wedge e_5^*$ and $\mgoth{b}= span (e_2,e_3,e_4,e_5 )$ with Lie bracket $[e_2,e_5]= e_3; [e_3,e_5]= -e_2.$
\end{description}

\noindent{\bf Semi-direct products}
\begin{description}
 \item [D4] $  [e_2,e_3] = e_1$, $[e_1,e_5] = (1+p)e_1$, $[e_2,e_5] = e_2$, $[e_3,e_5] = pe_3$, $[e_4,e_5] = -2(p+1)e_4$,
  $p\neq -1$, $\eta=e_1^*+e_4^*$.
  Here $\mgoth{g}$ is the semi-direct product $(\mgoth{h}_3\oplus
	     \mathbb Re_4)\rtimes \mathbb Re_5$ where $\mgoth{h}_3\oplus \mathbb Re_4$ is  as in  Section \ref{Heisenberg}.

 \item [D8] \ $[e_2,e_3]= e_1; [e_1,e_5]= 2e_1; [e_2,e_5]= e_2 +  e_3; [e_3,e_5]= e_3; [e_4,e_5]=  -4e_4$;
  $\eta=e^*_1+e^*_4$.  This is the semi-direct product $(\mgoth{h}_3\oplus \mathbb Re_4)\rtimes \mathbb Re_5$.
 
\item  [D10] \ $[e_2,e_3]= e_1; [e_1,e_5]= 2pe_1; [e_2,e_5]= pe_2+ e_3;
	     [e_3,e_5]=-e_2+ pe_3; [e_4,e_5]= -4p e_4$, $p\neq 0$;
	     $\eta=e^*_1 +e^*_4$. This is the semi-direct product $(\mgoth{h}_3\oplus \mathbb Re_4)\rtimes \mathbb Re_5$.

\item [D13] \ $ [e_2,e_3]= e_1; [e_1,e_5]= - \frac{1}{2}e_1; [e_2,e_5]=   - \frac{3}{2}e_2; [e_3,e_5]= e_3 + e_4; [e_4,e_5]=  e_4$; $\eta = e_1^* + e^*_4$; $p\neq 0$.  This is the semi-direct product $(\mgoth{h}_3\oplus \mathbb Re_4)\rtimes \mathbb Re_5$.

\item [D15] \ $[e_2,e_4] = e_1$, $[e_3,e_4] = e_2$, $[e_1,e_5] =\frac{2}{3}e_1$, $[e_2,e_5] =  -\frac{1}{3}e_2$, 
$  [e_3,e_5] =   - \frac{4}{3} e_3$,
$  [e_4,e_5] = e_4$,  $\eta=e^*_1+ e_3^*$. This is the  semi-direct product
$\mgoth{b}\rtimes \mathbb Re_5$ where $\mgoth{b}$ is the nilpotent
Lie algebra $\mgoth{b}=span(e_1,e_2,e_3, e_4)$ (Note that this is the 15th entry of the list in \cite{Diatta-Contact} with $p=-\frac{4}{3}.$)
\item [D18]
 $[e_1,e_4] = e_1$,  $[e_3,e_4] = - e_3$, $[e_2,e_5] = e_2$,  $[e_3,e_5] = -e_3$; $\eta=e^*_1+e^*_2+e^*_3$.
\item [D20] $[e_1,e_4]=-2e_1;$ $ [e_2,e_4]=e_2;$ $ [e_3,e_4]=e_3$; $[e_2,e_5]=-e_3;$ $ [e_3,e_5]=e_2.$

The last two Lie algebras above are the $2$-step solvable Lie algebra  $\mathcal R^3\rtimes \mathbb R^2$ where the Abelian subalgebra $\mathcal R^3=span (e_1,e_2,e_3)$ is the derived ideal and   $\mathcal R^2=span(e_4,e_5)$ is also Abelian.
\end{description}

Inspection of the list above yields the following corollary.

\begin{corollary}
Let $G$ be a five-dimensional simply-connected solvable contact Lie group.
Then $G$ is splittable.
\end{corollary}

See Appendix I for a list of descriptions of the nilradicals for each of these Lie groups.


\section{Five-dimensional contact Lie groups with uniform lattices}\label{chap:main results}

The following theorem indicates which of the simply-connected contact Lie groups in Theorem \ref{thm:Diatta1} have uniform lattices.

\begin{theorem}\label{thm:listing}Let $G$ be a five-dimensional connected and simply connected contact Lie group with a uniform lattice.  Then one of the following statements is true.
\begin{enumerate}

\item $G$ is the central extension of a solvable symplectic Lie group with a lattice that extends to $G.$  In particular, $G$ is one of the following groups:

\begin{enumerate}

\item $\heis^5=\R^4\times_{\omega_1}\R,$ where $\omega_1$ is the standard symplectic form on $\R^4,$

\item $\left(\heis^3\times \R\right)\times_{\omega_2} \R,$ where $\omega_2$ is the symplectic form on $\heis^3\times\R,$\ or

\item $B_j\times_{\omega_{j}} \R$\ ($j=3,\ 4,\ 5$), where $\omega_j$ is the symplectic form on $B_j=\R^3\rtimes_{F_j} \R$ with $F_j:\R\rightarrow Gl(3,\R)$ defined by the matrices
\begin{enumerate}
\item $F_3(t)=\left(
\begin{array}{ccc}1&-t&\frac{1}{ 2}t^2\\ 0&1&-t\\ 0&0&1
\end{array}
\right),$
\item $F_4(t)=\left(
\begin{array}{ccc}e^{-t}&0&0\\ 0&e^t&0\\ 0&0&1
\end{array}
\right),$
\item $F_5(t)=\left(
\begin{array}{ccc}cos(t)&-sin(t)&0\\ sin(t)&cos(t)&0\\ 0&0&1
\end{array}
\right).$
\end{enumerate}

\end{enumerate}

\item $G$ is a solvable semi-direct product and is one of the following groups:
\begin{enumerate}
\item $\R^3\rtimes_{b_1} \R^2,$ where $b_1:\R^2\rightarrow Gl(3,\R)$ is given  by $b_1(s,t)= \begin{pmatrix} e^{-s}&0&0
\\
 0  &e^{-t}& 0
\\
 0&0&e^{s+t}
\end{pmatrix}.$
\item $\R^3\rtimes_{b_2} \R^2,$ where $b_2:\R^2\rightarrow Gl(3,\R)$ is given  by
$b_2(s,t)= \begin{pmatrix} e^{2s}&0&0\\ 0 &e^{-s}\cos(t)&-e^{-s}\sin(t)\\ 0&e^{-s}\sin(t)&e^{-s}\cos(t)\\
\end{pmatrix}.$
\end{enumerate}
\end{enumerate}
\end{theorem}

\begin{corollary}
Let $X$ be a compact five-dimensional contact manifold uniformized by a five-dimensional contact Lie group $G.$  Then $G$ is solvable.
\end{corollary}


\section{Proof of Theorem \ref{thm:listing}}\label{chap:big proof}

We will prove Theorem \ref{thm:listing} by showing (1) that the Lie groups stated in the Theorem have lattices and (2) that the rest of the five-dimensional contact Lie groups given by Theorem \ref{thm:Diatta1} and the list in Subsection \ref{5-dimensional-unimod} do not.  As stated before, each solvable Lie group (or Lie algebra) will be referred to by its label in the list, e.g. {\bf D2}, {\bf D13}.  For ease of reading, we have relegated several technical results to appendices at the end of the paper.  In Appendix I (Section \ref{nillist}), the reader will find   a description of the nilradical of each solvable Lie algebra in the list in Subsection
\ref{5-dimensional-unimod} as well as matrix representations of $db$ and $\beta$ for the splitting $\mgoth{n}\rtimes_{\beta}T.$

\subsection{Positive cases}
The groups listed in Theorem \ref{thm:listing} are the simply connected Lie groups with Lie algebras {\bf D1}, {\bf D2}, {\bf D3}, {\bf D5}, {\bf D11}, {\bf D18}, and {\bf D20}, respectively.
Due to the variety of specific procedures used, we prove the existence of lattices in several
individual propositions. The overall methodology for the non-nilpotent cases is that utilized by Sawai in \cite{sawai} and by Sawai and Yamada in \cite{sy}, in which the existence of specific lattices is proven. 

\vskip .2in

\begin{proposition}
The Lie groups with Lie algebras {\bf D1}, {\bf D2} and {\bf D3} have lattices.
\end{proposition}

\vskip .2in

\noindent{\bf Proof:}\ \
 The Lie algebras {\bf D1}, {\bf D2} and {\bf D3} are all nilpotent.
 Recall that a nilpotent Lie group has a lattice if and only if its Lie
 algebra  also has a $\Q$-algebra, i.e., there is a basis on which all the coefficients of all of the bracket relations are in $\Q$ (See pp. 46-47 of \cite{ov1}).  The bases of the Lie algebras for these Lie groups as given in Appendix I (Section \ref{5-dimensional-unimod}) all satisfy this property.  Thus, groups with Lie algebras  {\bf D1}, {\bf D2} and {\bf D3} have lattices. \qed

\vskip .2in

\begin{proposition}
The Lie group with Lie algebra {\bf D5} has a lattice.
\end{proposition}
\noindent{\bf Proof:}\ \ Let $G$ be the simply-connected, connected Lie group with
Lie algebra {\bf D5}.  Then $G=(\heis^3 \times \R) \rtimes_b \R.$  For this proposition, we use
a representation of the multiplication on $\heis^3$ different from that described in the previous section.
Namely, for $\begin{pmatrix} x^\pr \\ y^\pr \\ z^\pr \\ \end{pmatrix},\ 
\begin{pmatrix} x \\ y \\ z \\ \end{pmatrix} \in \heis^3,$ let
$$ \begin{pmatrix} x^\pr \\ y^\pr \\ z^\pr \\ \end{pmatrix} \cdot 
\begin{pmatrix} x \\ y \\ z \\ \end{pmatrix} 
=\begin{pmatrix} x^\pr + \frac{1}{2}(y^\pr z - y z^\pr) +x \\ y^\pr +y\\ z^\pr +z\\ \end{pmatrix}.$$
That is, 
the $x$-coordinate of the product is the sum of the two $x$-coordinates plus one half of the signed area
of the parallelogram in the plane defined by the origin, $\begin{pmatrix} y^\pr \\ z^\pr \\ \end{pmatrix}$ and
$\begin{pmatrix} y \\ z \\ \end{pmatrix}$.  
Set $A\left(\begin{pmatrix} y^\pr \\ z^\pr \\ \end{pmatrix},\begin{pmatrix} y \\ z \\ \end{pmatrix}\right)= \frac{1}{2}(y^\pr z - y z^\pr).$
With this version of the Heisenberg group, multiplication on $G$ is given by
$$  
\begin{pmatrix} x^\pr \\ y^\pr \\ z^\pr \\ w^\pr\\ t^\pr \\ \end{pmatrix} \cdot
\begin{pmatrix} x \\ y  \\ z  \\ w \\ t  \\ \end{pmatrix} =
\begin{pmatrix} x^\pr - \frac{1}{2}( e^{t^\pr}z^\pr  y -  e^{-t^\pr}y^\pr   z) -t^\pr w +x\\ 
y^\pr+e^{-t^\pr} y \\ z^\pr +e^{t^\pr} z\\ w^\pr+w\\ t^\pr +t \\ \end{pmatrix}
$$

Let $n$ be a positive integer greater than 2.
Let $t_0$ be a positive real number such that $(e^{t_0})^2 - n e^{t_0}+1=0$.  
Set $v_1=\begin{pmatrix} 1 \\ 1\\ \end{pmatrix},\ 
v_2=\begin{pmatrix} e^{-t_0} \\ e^{t_0}\\ \end{pmatrix}\in \R^2$ and 
$x_0= A(v_1, v_2).$
Note that the lattice of $\R^2$ generated by the vectors $v_1$ and $v_2$
is preserved by the linear transformation $\begin{pmatrix} e^{-t_0}& 0\\ 0& e^{t_0} \\ \end{pmatrix}.$

Let $\Gamma\in G$ be the discrete subset given by
$$\Gamma=
\Z\begin{pmatrix} x_0 \\ 0 \\ 0  \\ 0 \\ 0  \\ \end{pmatrix}+
\Z\begin{pmatrix} 0 \\ v_1  \\ 0 \\ 0  \\ \end{pmatrix} +
\Z\begin{pmatrix} 0 \\ v_2  \\ 0 \\ 0  \\ \end{pmatrix}+
\Z\begin{pmatrix} 0 \\ 0 \\ 0  \\ \frac{x_0}{t_0} \\ 0  \\ \end{pmatrix}+
\Z\begin{pmatrix} 0 \\ 0 \\ 0  \\ 0 \\ t_0  \\ \end{pmatrix}.$$
By definition of the quantities and vectors involved, we can see that $\Gamma$ is actually a {\it subgroup} of $G$ (c.f. Sawai \cite{sawai} and Sawai and Yamada \cite{sy}).
Thus, $G$ has a lattice. \qed

\begin{proposition}
The Lie group with Lie algebra {\bf D11} has a lattice.
\end{proposition}
\noindent{\bf Proof:}\ \ Let $G$ be the simply-connected, connected Lie group with Lie algebra {\bf D11}.
The group structure of $G$ is given by
$$\begin{pmatrix} x_1 \\ y_1\\ z_1\\ w_1\\ \theta_1\\
\end{pmatrix} \cdot
\begin{pmatrix} x_2 \\ y_2\\ z_2\\ w_2\\ \theta_1\\
\end{pmatrix} =
\begin{pmatrix} x_1 +x_2 \pm \theta_1 w_2 +y_1((\sin \theta_1)y_2 + (\cos \theta_1)z_2)           \\ 
y_1+ (\cos \theta_1)y_2 - (\sin \theta_1)z_2\\ 
z_1+ (\sin \theta_1)y_2 + (\cos \theta_1)z_2) \\ 
w_1+w_2\\ 
\theta_1+\theta_2\\
\end{pmatrix}. $$
Then the subgroup of $G$ generated by the elements
$$\left\{    
\begin{pmatrix} \pi \\ 0\\ 0\\ 0\\ 0\\
\end{pmatrix},
\begin{pmatrix} 0\\ \sqrt{\pi} \\ 0\\ 0\\ 0\\
\end{pmatrix},
\begin{pmatrix} 0\\ 0\\ \sqrt{\pi} \\ 0\\ 0\\
\end{pmatrix},
\begin{pmatrix} 0\\ 0 \\ 0\\ 1\\ 0\\
\end{pmatrix},
\begin{pmatrix} 0\\ 0 \\ 0\\ 0\\ \pi \\
\end{pmatrix}     
\right\}$$
is discrete in $G$ and hence a lattice. \qed

\begin{proposition}
The Lie group with Lie algebra {\bf D18} has a lattice.
\end{proposition}
\noindent{\bf Proof:}\ \ 
Let $G=\R^3 \rtimes_{b_1} \R^2$ be the Lie group corresponding to Lie algebra {\bf D18}.
In order to show that $G$ has a lattice, we need to produce a basis of $\R^3$, $\{ v_1, v_2, v_3\}$, and a basis of $\R^2$, $\{p_1, p_2\}$ such that
$$\Gamma = \left< v_1, v_2, v_3 \right>_\Z \rtimes_{b_1} \left< p_1, p_2\right>_\Z$$
is a subgroup of $G.$
We do this by setting 
$$T_1=\left(
\begin{array}{ccc}0&0&1\\ 1&0&-5\\ 0&1&6\\
 \end{array}
\right),\ \ T_2=\left(
\begin{array}{ccc}-4&-4&-3\\ 21&16&11\\ -4&-3&-2\\
\end{array}
\right).$$
It is easily verified that $T_1 \circ T_2 = T_2\circ T_1.$

\noindent Furthermore, 
 the characteristic polynomials of $T_1$ and $T_2$ are given respectively as
\begin{eqnarray*}
f_1(X)&=&X^3-6X^2+5X-1, \\
f_2(X)&=&X^3-10X^2+17X-1, \\
\end{eqnarray*}
each of which have three distinct roots.  
The roots of $f_1$ are
\begin{eqnarray*}
\alpha_1 &=& 0.30797853...\\
\beta_1 &=& 0.64310413...\\
\gamma_1 &=&  5.0489173...,\\
\end{eqnarray*}
and the roots of $f_2$ are
\begin{eqnarray*}
\alpha_2 &=& 2.088146...\\
\beta_2 &=& 7.8508551...\\
\gamma_2 &=&  0.06099892...,\\
\end{eqnarray*}
Thus, $T_1$ and $T_2$ are simultaneously diagonalizable.
In fact, there is 
a $\Phi\in Gl(3,\R)$ such that, for $j=1,\ 2,$
$$\Phi T_j \Phi^{-1} = \left(\begin{array}{ccc}\alpha_j&0&0\\0 &\beta_j&0\\ 0 &0&\gamma_j\\
\end{array}  \right).$$

Set 
$$
p_1 = \left( \begin{array}{c} \ln \alpha_ 1\\  \ln \beta_ 1\\ \end{array} \right) 
= \left( \begin{array}{c} -1.7772... \\ -0.441449... \\ \end{array}\right),\ \ \ 
p_2 = \left( \begin{array}{c} \ln \alpha_2 \\  \ln \beta_2 \\ \end{array} \right) 
= \left(\begin{array}{c} 2.06062...\\  0.736277...\\  \end{array} \right),\ \ 
$$
Note that the slope of $p_1$ as a vector is approximately $2.66786$ and that of $p_2$ is approximately $2.79871.$  Thus, $p_1$ and $p_2$ are linearly independent vectors in $\R^2.$
By definition of $p_1$ and $p_2$, we have
$$b(p_j)= \left(\begin{array}{ccc}\alpha_j&0&0\\0 &\beta_j&0\\ 0 &0&\gamma_j\\
\end{array}  \right)$$
for each $j=1,2$  so that $b(p_j) \circ \Phi = \Phi \circ T_j$ for each $j.$
In particular, if we set $v_k = \Phi\left(\epsilon_k \right)$ for $k=1,\ 2,\ 3,$
where $\{\epsilon_1, \epsilon_2, \epsilon_3\}$ is the standard basis of $\R^3,$
then $b(p_j)(v_k) \in \left< v_1, v_2, v_3\right>_Z$ for each $j=1,\ 2$ and $k=1, 2, 3.$
Thus, $\Gamma = \left< v_1, v_2, v_3 \right>_\Z \rtimes_{b_1} \left< p_1, p_2\right>_\Z$
is a discrete, uniform subgroup of $G.$  This proves the claim.

\begin{proposition}
The Lie group with Lie algebra {\bf D20} has a lattice.
\end{proposition}
\noindent{\bf Proof:}\ \ 
Let $G=\R^3 \rtimes_{b_1} \R^2$ be the Lie group corresponding to Lie algebra {\bf D20}.
As with the previous claim, we need to 
produce a basis of $\R^3$, $\{ v_1, v_2, v_3\}$, and a basis of $\R^2$, $\{p_1, p_2\}$ such that
$$\Gamma = \left< v_1, v_2, v_3 \right>_\Z \rtimes_{b_2} \left< p_1, p_2\right>_\Z$$
is a subgroup of $G.$

Here, we do this by setting
$$U_1=\left( \begin{array}{ccc} 0 & 0 & 1\\ 1&0&-2\\ 0&1&3\\
\end{array}\right),\ \ \ U_2=\left( \begin{array}{ccc} 0 & 1 & 1\\ -2&-2&-1\\ 1&1&1\\
\end{array}\right).$$
Then $U_1 U_2= U_2 U1$.
Also,  the characteristic polynomials of $U_1$ and $U_2$ are given respectively as
\begin{eqnarray*}
f_1(X)&=&X^3-3X^2+2X-1, \\
f_2(X)&=&X^3+X^2-1, \\
\end{eqnarray*}
each of which have exactly one real root -- $f_1(2.3247...)=0=f_2(0.7549...)$ -- and two
complex roots.

Thus, there is 
a $\Psi\in Gl(3,\R)$ such that, for $j=1,\ 2,$
$$
\Psi U_j \Psi^{-1}=\left(\begin{array}{ccc}{\alpha_j}^2&0&0
\\
0&{\alpha_j}^{-1}cos(\beta_j)&-{\alpha_j}^{-1}sin(\beta_j)
\\  0&{\alpha_j}^{-1}sin(\beta_j)&{\alpha_j}^{-1}cos(\beta_j)
\\
\end{array}  \right),
$$
where $\alpha_1\approx \sqrt{2.3247},$ $\beta_1\approx 1.0300,$
$\alpha_2\approx \sqrt{0.7549}$ and $\beta_2\approx 2.4378.$

Set 
$$
p_1 = \left( \begin{array}{c} \ln \alpha_ 1\\  \beta_ 1\\ \end{array} \right) 
= \left( \begin{array}{c} 0.4217... \\ 1.0300... \\ \end{array}\right),\ \ \ 
p_2 = \left( \begin{array}{c} \ln \alpha_2 \\  \beta_2 \\ \end{array} \right) 
= \left(\begin{array}{c} -0.1405...\\  2.4378...\\  \end{array} \right),\ \ 
$$
Thus, $p_1$ and $p_2$ are linearly independent vectors in $\R^2.$
By definition of $p_1$ and $p_2$, we have
$$b(p_j)= \left(\begin{array}{ccc}{\alpha_j}^2&0&0
\\
0&{\alpha_j}^{-1}cos(\beta_j)&-{\alpha_j}^{-1}sin(\beta_j)
\\  0&{\alpha_j}^{-1}sin(\beta_j)&{\alpha_j}^{-1}cos(\beta_j)
\\
\end{array}  \right)$$
for each $j=1,2$  so that $b(p_j) \circ \Psi = \Psi \circ U_j$ for each $j.$
In particular, if we set $v_k = \Psi\left(\epsilon_k \right)$ for $k=1,\ 2,\ 3,$
where $\{\epsilon_1, \epsilon_2, \epsilon_3\}$ is the standard basis of $\R^3,$
then $b(p_j)(v_k) \in \left< v_1, v_2, v_3\right>_Z$ for each $j=1,\ 2$ and $k=1, 2, 3.$
Thus, $\Gamma = \left< v_1, v_2, v_3 \right>_\Z \rtimes_{b_1} \left< p_1, p_2\right>_\Z$
is a discrete, uniform subgroup of $G.$  This proves the claim.

\subsection{Negative cases}

\subsubsection{Solvable}

We now show that the rest of the solvable contact Lie groups of five dimensions do not have lattices.  Appendix I lists the corresponding homomorphism $\beta:T\rightarrow der(\mgoth{n})$ and $db:T\rightarrow Aut(\mgoth{n})$ of each group.  
  There are only two classes of such Lie groups remaining from the list in Subsection \ref{5-dimensional-unimod}, namely, those whose nilradical is $\heis^3\times \R$
  ({\bf D4,\ D5,\ D8,\ D10,\ D11,\ D13}) and one whose nilradical is a
  semidirect product $\R^3\rtimes_f \R$ ({\bf D15}).

\medskip

\noindent

\begin{proposition}None of the Lie groups corresponding to the Lie algebras
{\bf D4,\ D8,\ D10,\ D13} have lattices.
\end{proposition}

\noindent {\bf Proof:}\ \ 
Let $\mgoth{g}$ be one of the Lie algebras {\bf D4,\ D8,\ D10,\ D13}.
Let $G$ be the simply-connected Lie group with Lie algebra $\mgoth{g}.$
Then we have the short exact sequence
$$0\rightarrow N \rightarrow G \rightarrow \R \rightarrow 0,$$
where $N=\heis^3 \times \R$ and the sequence splits ($G=N \rtimes_b \R$).

Suppose $G$ has a lattice $\Gamma.$ 
By Theorem \ref{thm:mostow}, $N\cap\Gamma$ is a lattice of $N$; and, 
by Theorem \ref{thm:wang}, $\Gamma$ is isomorphic to a group of the form
$(N\cap\Gamma) \rtimes_b \left< t_0 \right>_\Z$ for some $t_0\in \R^+.$
Recall also that a lattice of $N$ necessarily contains a subgroup given as the
span over $\Z$ of 
$g_0=\begin{pmatrix} x_0\\ 0\\ 0\\ 0\\ \end{pmatrix}$ for some positive $x_0\in \R.$
Thus,
$g_1=b(t_0)(g_0)  $ 
and ${g_1}^{-1}=b(-t_0)(g_0)  $ 
are both elements of $N\cap\Gamma$.

We now show the non-existence of lattices in $G$ by looking at each of the possible cases. 

\begin{enumerate}[\ \ \hskip 25pt \ \ \ A.]
\item If $\mgoth{g}$ is {\bf D4}, 
then 
$$g_1 = \begin{pmatrix} e^{-(p+1)t_0}x_0\\ 0\\ 0\\ 0\\ \end{pmatrix}\ {\rm and}\ 
g_2 = \begin{pmatrix}e^{(p+1)t_0} x_0\\ 0\\ 0\\ 0\\ \end{pmatrix}.$$
Since $g_1$ and $g_2$ are integer multiples of $g_0$, we have either $p=-1$ or $t_0=0.$ 
Both of these possible cases are contradictions.

\item If $\mgoth{g}$ is {\bf D8}, 
then 
$$g_1 = \begin{pmatrix} e^{-2t_0}x_0\\ 0\\ 0\\ 0\\ \end{pmatrix}\ {\rm and}\ 
g_2 = \begin{pmatrix}e^{2t_0} x_0\\ 0\\ 0\\ 0\\ \end{pmatrix}.$$
Since $g_1$ and $g_2$ are integer multiples of $g_0$, we have $t_0=0.$ 
This is a contradiction.

\item If $\mgoth{g}$ is {\bf D10}, 
then 
$$g_1 = \begin{pmatrix} e^{-2pt_0}x_0\\ 0\\ 0\\ 0\\ \end{pmatrix}\ {\rm and}\ 
g_2 = \begin{pmatrix}e^{2pt_0} x_0\\ 0\\ 0\\ 0\\ \end{pmatrix}.$$
Since $g_1$ and $g_2$ are integer multiples of $g_0$, we have either $p=0$ or $t_0=0.$ 
Both of these possible cases are contradictions.

\item 
Finally, if $\mgoth{g}$ is {\bf D13}, 
then 
$$g_1 = \begin{pmatrix} e^{\frac{1}{2} t_0}x_0\\ 0\\ 0\\ 0\\ \end{pmatrix}\ {\rm and}\ 
g_2 = \begin{pmatrix}e^{-\frac{1}{2} t_0} x_0\\ 0\\ 0\\ 0\\ \end{pmatrix}.$$
Since $g_1$ and $g_2$ are integer multiples of $g_0$, we have $t_0=0.$ 
This is a contradiction.
\end{enumerate}
This proves the proposition.\qed

\begin{proposition}
A Lie group with Lie algebra {\bf D15} does not have a lattice.
\end{proposition}
\noindent{\bf Proof:}\ \ The simply connected Lie group with Lie algebra {\bf D15} is given by
$G=(\R^3 \rtimes_f \R)\rtimes_b \R,$
where $f: \R \rightarrow \R^3$ is given by
$$f(w)=\begin{pmatrix} 1 & -w & \frac{1}{2} w^2 \\
&1&-w\\ &&1 \end{pmatrix} $$
and $b: \R \rightarrow \R^3 \rtimes_f \R$ is given by
$$b(t)=\begin{pmatrix} e^{-\frac{2}{3}t} & & & \\
&e^{\frac{1}{3}t} & & \\ & & e^{\frac{4}{3}t}& \\ & & & e^{-t}\\ \end{pmatrix}. $$

The nilradical of group {\bf D15} is of the subgroup $N= (\R^3 \rtimes_f \R)\rtimes_b (0) $
with multiplication given by
$$\begin{pmatrix}x \\ y\\ z\\ w\\    \end{pmatrix} \cdot
\begin{pmatrix}x^\pr \\ y^\pr\\ z^\pr\\ w^\pr\\    \end{pmatrix}    =
\begin{pmatrix}x+x^\pr -wy^\pr +\frac{1}{2}w^2z^\pr \\y+ y^\pr +wz^\pr \\ 
z+z^\pr\\ w+w^\pr\\    \end{pmatrix}. 
$$
so that 
$$ \begin{pmatrix} x\\ y\\ z\\ w\\ \end{pmatrix}^{-1} = \begin{pmatrix} -x-yw-\frac{1}{2}zw^2\\ -y-zw\\ -z\\ -w\\ \end{pmatrix}.  $$
And, for $g_j=\begin{pmatrix} x_j\\ y_j\\ z_j\\ w_j\\ \end{pmatrix}\in N,$ $j=1,2,$
$$g_1 g_2 {g_1}^{-1} {g_2}^{-1}
=\begin{pmatrix}y_1 w_2-y_2w_1 -\frac{1}{2} z_1 {w_2}^2 +\frac{1}{2} z_2 {w_1}^2\\ 
z_1w_2-z_2w_1\\ 0\\ 0\\    \end{pmatrix}.$$
Let $N_1=[N,N]$ and $N_2=[[N,N],N]$.
Then $N_2=\R\times (0,0,0),$ and
$[N_2,N]=(0),$ that is, $N$ is 3-step nilpotent.

Suppose $G$ has lattice $\Gamma.$  By Theorem \ref{thm:mostow}, 
$\Gamma_0 = \Gamma \cap N$ is a lattice of $N$. By Theorem \ref{handy}, 
$\Gamma_1 = \Gamma \cap N_1$ is a lattice
of $N_1,$ and $\Gamma_2 = \Gamma \cap N_2$ is a lattice
of $N_2\cong \R$.  Let $x_0$ be the unique positive real number such that 
$$\Gamma_2=\left\{\begin{pmatrix} n\ x_0 \\ 0 \\ 0\\ 0\\ \end{pmatrix}: n\in \Z\right\}.$$

Furthermore, by Theorem \ref{thm:wang}, $\Gamma$ is isomorphic to a group
$\tilde\Gamma$ satisfying the short, exact sequence
$$0\rightarrow  \Gamma_0 \rightarrow  \tilde\Gamma \rightarrow \Gamma_0 \setminus  
\tilde\Gamma  \rightarrow 0$$
induced from the short exact sequence
$$0\rightarrow N \rightarrow G  \rightarrow   N\setminus G\rightarrow 0.$$
Since $N\setminus G\cong \R,$ there is a $t_0\in \R^+$ such that
$b(\pm t_0)(\Gamma_0) \subset \Gamma_0.$
Since $ b(t_0)$ and $ b(-t_0)$ are non-trivial group homomorphisms from $\Gamma_0$ to itself, 
they must preserve the central series of $\Gamma_0,$ that is,
$$ b(\pm t_0)(\Gamma_1) \subset \Gamma_1\ {\rm and}\  
b(\pm t_0)(\Gamma_2) \subset \Gamma_2.$$
Thus, both
$$b(t_0)\begin{pmatrix} x_0\\ 0\\ 0\\ 0\\ \end{pmatrix} =
\begin{pmatrix}e^{-\frac{2}{3}t_0} x_0\\ 0\\ 0\\ 0\\ \end{pmatrix},\ 
b(-t_0)\begin{pmatrix} x_0\\ 0\\ 0\\ 0\\ \end{pmatrix} =
\begin{pmatrix}e^{\frac{2}{3}t_0} x_0\\ 0\\ 0\\ 0\\ \end{pmatrix}
\in \left\{\begin{pmatrix} n\ x_0 \\ 0 \\ 0\\ 0\\ \end{pmatrix}: n\in \Z\right\},
$$
implying that $e^{\pm\frac{2}{3}t_0}\in \Z.$  So, $t_0=0$, which is a contradiction, and thus no lattice exists on $G.$
\qed

\subsubsection{Non-solvable: the general case of $\R^n \rtimes Sl(n,\R).$} \label{specialaffine}

According to Theorem \ref{thm:Diatta1}, the
 only unimodular non-solvable contact Lie group of dimension five is the
 group 
 $\R^2 \rtimes Sl(2,\R)$ of special affine transformations of the
 plane. 
We obtain the following more general result stating the nonexistence of uniform
 lattices in  $\R^n \rtimes Sl(n,\R)$, for every $n\ge 2$. See Theorem
 \ref{nolatticesln}. Let us recall that  $\R^n \rtimes Sl(n,\R)$ is a
 contact Lie group \cite{Diatta-Contact}. We can exhibit a contact form $\eta$
 on the Lie algebra  $\R^n \rtimes sl(n,\R)$ of  $\R^n \rtimes Sl(n,\R)$
 by looking at it as the subalgebra 
 $\R^n \rtimes sl(n,\R)=\left\{ \begin{pmatrix} A & v\\ 0&0\end{pmatrix}, ~
 \text{ where $A\in sl(n,\R)$  and $v\in\R^n$} \right\}$ of the Lie
 algebra $\mathcal Gl(n+1,\R)$ of $(n+1)\times (n+1)$ real matrices.  The
 $(n+1)\times (n+1)$ matrices $e_{i,j}$ all of whose entries are zero except the ${ij}$-th one  which is
 equal to $1$, form a basis of  $\mathcal Gl(n+1,\R)$. Let us denote by
 ($e_{i,j}^*$) the corresponding dual basis.
Then, $\eta:=\displaystyle \sum_{i=1}^n e_{i,i+1}^*$  is a contact form
 on $\R^n \rtimes sl(n,\R)$, with Reeb vector $\xi:=\displaystyle
 \frac{1}{n}\sum_{i=1}^n e_{i,i+1}.$ Now, we have the following.

\begin{theorem}\label{nolatticesln}
The group  $\R^n \rtimes Sl(n,\R)$, of special affine transformations of
 $\R^n$, has no uniform lattice, for every $n\ge 2$.
\end{theorem}

\noindent{\bf Proof:}~~
 Let $G:= \R^n \rtimes Sl(n,\R)$ and suppose $\Gamma$ is a lattice in
 $G$.  The radical of $G$ is the subgroup
 $\R^n\times \{I\}.$  Then $\Gamma^\pr=\Gamma \cap \R^n\times \{I\}$ is
 a lattice of $\R^n\times \{I\}$ (Corollary 1.8 on p. 107 of
 \cite{ov2}). Let $v_1, \ldots,  v_n\in\R^n$ be such that $(v_1,I)\ldots
 (v_n,I)$  generate $\Gamma^\pr.$  Let $A\in Sl(n,\R)$ and $w\in \R^n$
 such that $(w,A)\in\Gamma.$  Then, for $j=1, \ldots, n,$
\begin{eqnarray*} (w,A)(v_j,I)(w,A)^{-1}=(Av_j+w,A)(-A^{-1}w,A^{-1})=(Av_j,I)\in\Gamma.
\end{eqnarray*} Hence the set $M_\Gamma$ given by $$M_\Gamma=\{A\in
 Sl(n,\R): (w,A)\in\Gamma\ {\rm for\ some\ }w\in \R^n\}$$ preserves the
 lattice $\Gamma^\pr$ on $\R^n.$  In particular, by the change of basis
 $v_j\mapsto e_j$ for $j=1,\ldots, n,$ we can assume that
 $M_\Gamma\subset Sl(n,\Z).$  Now, it is known that $Sl(n,\Z)$ is a
 lattice of $Sl(n,\R)$ but not a uniform lattice (e.g., see pp. 229-231 of \cite{beardon} for the case where $n=2$). In other words, there is a sequence $\{\gamma_j\}\subset Sl(n,\R)$ such that its projection $Sl(n,\Z)\setminus Sl(n,\R)$ has no convergent subsequences.  Thus, its projection in $M_\Gamma\setminus Sl(n,\R)$ also has no convergent subsequences, which means that the sequence $\{[0,\gamma_j]\}\subset \Gamma\setminus \R^n \rtimes Sl(n,\R)$ has no convergent subsequences.  Therefore, $\Gamma\setminus\R^n \rtimes Sl(n,\R)$ is not compact.  Since $\Gamma$ was assumed to be an arbitrary lattice of $\R^n \rtimes Sl(n,\R)$, $\R^n \rtimes Sl(n,\R)$ has no uniform lattices.\qed

\section{Appendix I: List of nilradicals of the unimodular contact Lie algebras of dimension $5$ }\label{nillist}

The following is a list of all of the unimodular Lie algebras among
those in  the first author's list of solvable contact Lie groups in five
dimensions from \cite{Diatta-Contact}. Their Lie brackets, in a basis $ (e_1, e_2, e_3, e_4, e_5)$,  are given in Section  \ref{5-dimensional-unimod}. 
Each of the corresponding Lie groups will be of the form $N\rtimes_b T$,
where $N$ is the nilradical, $T$ is an Abelian group and $b:T\rightarrow
Aut(N)$ a homomorphism. For each of these, the Lie algebra $\mgoth{n}$
of the nilradical $N$ of the simply-connected Lie group corresponding to
each Lie algebra is provided as well as the Abelian group $T$.  The
transformations $\beta$ and $db$ are matrix representations (with
respect to the given basis of $\mgoth{n}$) of the corresponding
homomorphisms $\beta:T\rightarrow der(\mgoth{n})$ and $db(x)=exp(\beta(x)):T\rightarrow
Aut(\mgoth{n})$ (for $x\in T$) induced from the semidirect product
$N\rtimes_b T.$ 

 \begin{description}
\item [D1] 
$\mgoth{n}=<e_1, \dots, e_5>={\mathcal H}_5,$  $T=(0),$

\item [D2]
  $\mgoth{n}=(<e_1, e_3, e_4>\oplus<e_2>)+_{df}<e_5>=({\mathcal H}_3\oplus
	    \R)+_{df}\R,$ $T=(0)$  where
 $$
ad(e_5)=\left(\begin{array}{cccc} 0&0&0&-1\\ 0&0&0&0\\ 0&0&0&0\\ 0&-1&0&0\\ \end{array}\right),\ df(te_5)=\left(\begin{array}{cccc} 1&\frac{1}{2}t^2&0&-t\\ 0&1&0&0\\ 0&0&1&0\\ 0&-t&0&1\\ \end{array}\right) .
$$
\item [D3]
$
\mgoth{n}=(<e_1, e_3, e_4>\oplus<e_2>)+_{df}<e_5>=({\mathcal H}_3\oplus
	    \R)+_{df}\R,$  $T=(0)$ {\rm where} 
 $$ad(e_5)=\left(\begin{array}{cccc} 0&0&0&-1\\ 0&0&-1&0\\ 0&0&0&0\\ 0&-1&0&0\\ \end{array}\right),\ df(te_5)=\left(\begin{array}{cccc} 1&\frac{1}{2}t^2&-\frac{1}{6}t^3&-t\\ 0&1&-t&0\\ 0&0&1&0\\ 0&-t&0&1\\
\end{array}\right) .
$$
\item [D4]
$
\mgoth{n}=<e_1, e_2, e_3>\oplus< e_4>={\mathcal H}_3\oplus \R,$ $\ T=\R e_5, $
 $$\beta(e_5)=\left(\begin{array}{cccc} -(p+1)&0&0&0\\ 0&-1&0&0\\ 0&0&-p&0\\ 0&0&0&2(p+1)\\ \end{array}\right),\ db(te_5)=\left(\begin{array}{cccc} e^{-(p+1)t}&0&0&0\\ 0&e^{-t}&0&0\\ 0&0&e^{-pt}&0\\ 0&0&0&e^{2(p+1)t}\\ \end{array}\right).
$$
\item [D5]
$
\mgoth{n}=<e_1, e_2, e_3>\oplus< e_4>={\mathcal H}_3\oplus \R,\ T=\R e_5,
$
$$
\beta(e_5)=\left(\begin{array}{cccc} 0&0&0&-1\\ 0&-1&0&0\\ 0&0&1&0\\ 0&0&0&0\\ \end{array}\right),\ db(te_5)=\left(\begin{array}{cccc} 1&0&0&-t\\ 0&e^{-t}&0&0\\ 0&0&e^{t}&0\\ 0&0&0&1\\ \end{array}\right).
$$
\item [D8]
$
\mgoth{n}=<e_1, e_2, e_3>\oplus< e_4>={\mathcal H}_3\oplus \R,$ $ T=\R e_5,
$
 $$
\beta (e_5)=\left(\begin{array}{cccc} -2&0&0&0\\ 0&-1&0&0\\ 0&-1&-1&0\\ 0&0&0&4\\ \end{array}\right),\ db(te_5)=\left(\begin{array}{cccc}e^{-2t}&0&0&0\\ 0&e^{-t}&0&0\\ 0&-te^{-t}&e^{-t}&0\\ 0&0&0&e^{4t}\\ \end{array}\right) .
$$
\item  [D10] $
\mgoth{n}=<e_1, e_2, e_3>\oplus< e_4>={\mathcal H}_3\oplus \R, $ $ \ T=\R e_5,
$
$$
\beta(e_5)=\left(\begin{array}{cccc} -2p&0&0&0\\ 0&-p&1&0\\ 0&-1&-p&0\\ 0&0&0&4p\\
\end{array}\right),\ db(te_5)=
\left(\begin{array}{cccc}e^{-2pt}&0&0&0\\ 0&e^{-pt}\cos(-t)&-e^{-pt}\sin(-t)&0\\ 0&e^{-pt}\sin(-t)&e^{-pt}\cos(-t)&0\\  0&0&0&e^{4pt}\\
 \end{array}\right) .
$$
\item [D11] 
$
\mgoth{n}=<e_1, e_2, e_3>\oplus< e_4>={\mathcal H}_3\oplus \R,$ $\ T=\R e_5,
$
$$
\beta(e_5)=\left(\begin{array}{cccc} 0&0&0&\pm 1\\ 0&0&1&0\\ 0&-1&0&0\\ 0&0&0&0\\ \end{array}\right),\ db(te_5)= \left(\begin{array}{cccc}1&0&0&\pm t\\  0&cos(t)&-sin(t)&0 \\ 0&sin(t)&cos(t)&0 \\  0&0&0&1  \end{array}  \right).
$$
\item [D13] 
$
\mgoth{n}=<e_1, e_2, e_3>\oplus< e_4>={\mathcal H}_3\oplus \R,$ $T=\R e_5,
$
$$
\beta(e_5)=\left(\begin{array}{cccc} \frac{1}{2}&0&0&0\\ 0&\frac{3}{2}&0&0\\ 0&0&-1&0\\ 0&0&-1&-1
\\
\end{array}\right),\ db(te_5)=\left(\begin{array}{cccc}e^{\frac{1}{ 2}t}&0&0&0\\ 0&e^{\frac{3}{2}t}&0&0\\ 0&0&e^{-t}&0\\ 0&0&-te^{-t}&e^{-t}
\\
\end{array}\right) .
$$
\item [D15]
$
\mgoth{n}=<e_1, \dots, e_4>=<e_1,e_2,e_3>+_{f_*}<e_4>,$  ${\rm with}\ f_*(e_4)=\left(\begin{array}{ccc}0&-1&0\\ 0&0&-1\\ 0&0&0\\
\end{array}   \right),$ $ T=\R e_5,
$
$$
\beta (e_5)=\left(\begin{array}{cccc} -\frac{2}{3}&0&0&0\\ 0&\frac{1}{3}&0&0\\ 0&0&\frac{4}{3}&0\\ 0&0&0&-1\\ \end{array}\right),\ db(te_5)= \left(\begin{array}{cccc} e^{-\frac{2}{3}t}&0&0&0\\ 0&e^{\frac{1}{3}t}&0&0\\ 0&0&e^{\frac{4}{3}t}&0\\ 0&0&0&e^{-t}\\ \end{array}\right).
$$
\item [D18] $\mgoth{n}=<e_1, e_2, e_3 >$, with
$$
\beta (se_4+te_5)=\left( \begin{array}{ccc} -s&0&0\\ 0&-t&0\\ 0&0&s+t\\ \end{array}\right),\ db(se_4+te_5)=\left( \begin{array}{ccc} e^{-s}&0&0\\ 0&e^{-t}&0\\ 0&0&e^{s+t}\\ \end{array}\right).
$$
 \item [D20]  $\mgoth{n}=<e_1, e_2, e_3 >$, with
$$
\beta (se_4+te_5)=\left( \begin{array}{ccc} 2s&0&0\\ 0&-s&-t\\ 0&t&-s\\ \end{array}\right),\ db(se_4+te_5)=\left( \begin{array}{ccc} e^{2s}&0&0\\ 0&e^{-s}cos(t)&-e^{-s}sin(t)\\ 0&e^{-s}sin(t)&e^{-s}cos(t)\\ \end{array}\right)
$$
\end{description}

\medskip\noindent
{\bf Aknowledgement.}
The authors thank the referees of this paper for their attention to detail and excellent suggestions.  The second author thanks Dr. Patrick Chen at John Carroll University for his valuable insight on lattices and semi-direct products.


\addcontentsline{toc}{chapter}{References} \begin{thebibliography}{10}

  \bibitem {aus1} \textsc{L. Auslander}, An exposition of the structure of solvmanifolds I,\ Bull. Amer. Math. Soc., {\bf 79} (1973), 227-261.

\bibitem {aus2} \textsc{L. Auslander}, An exposition of the structure of solvmanifolds II,\ Bull. Amer. Math. Soc., {\bf 79} (1973), 262-285.

  \bibitem{barden} \textsc{D. Barden}, Simply connected five-manifolds.  Ann. of Math. (2), {\bf 82}  (1965), 365-385.

\bibitem{blz}\textsc{P. Basarab-Horwath and V. Lahno and R. Zhdanov}, The structure of Lie algebras and the classification problem for partial differential equations, Acta Appl. Math. {\bf 69} (2001), 43-94.

\bibitem{beardon} \textsc{A. Beardon}, The Geometry of Discrete Groups
 Graduate Texts in Mathematics {\bf 91}, Springer-Verlag (1983)

\bibitem{bg2} \textsc{ C. Benson and Gordon}, K\"ahler structures on compact solvmanifolds, Proc. American Math. Soc. {\bf 108} (1990), 971-980. 


\bibitem{blair} \textsc{D. Blair}, The Riemannian Geometry of Contact and Symplectic Manifolds, Birkh\"auser, New York, 2001.

\bibitem{borel}\textsc{ A. Borel}, Introduction aux groupes arithm\'ethiques, Hermann, Paris, 1969.

\bibitem{Diatta-Contact}  \textsc{A. Diatta}, Left Invariant Contact Structures on Lie Groups. Differential Geom. Appl. {\bf 26} (2008), 544-552.

\bibitem{Diatta-Contact-Riem}  \textsc{A. Diatta}, Riemannian Geometry on Contact Lie Groups. Geom. Dedicata {\bf 133} (2008), 83-94.

\bibitem{diatta-medina} \textsc{A. Diatta and A. Medina}, Classical Yang-Baxter Equation and Left Invariant Affine Geometry on Lie Groups. Manuscripta Math. {\bf 114} (2004), 477-486.

\bibitem{Geiges95} \textsc{H. Geiges}, Examples of symplectic $4$-manifolds with disconnected boundary of contact type. Bull. London Math. Soc. {\bf 27} (1995), 278-280

\bibitem{Geiges94} \textsc{H. Geiges},  Symplectic manifolds with disconnected boundary of contact type. Internat. Math. Res. Notices {\bf 1994} (1994), 23-30.

\bibitem{gw} \textsc{C. Gordon and E. Wilson},  The spectrum of the Laplacian on Riemannian Heisenberg manifolds , Michigan Math. J., {\bf 33} (1986), 253-271.



\bibitem{Harshavardhan} \textsc{R. Harshavardhan}, Geometric Structures of Lie Type on 5-Manifolds, Ph.D. Thesis Cambridge University, U.K. (1995-6).

\bibitem{mal}\textsc{ A. Mal'tsev}, Solvable Lie algebras, Amer. Math. Soc. Translation {\bf 1950}, (1950). no. 27.


\bibitem{mcduff}\textsc{D. McDuff},  Symplectic manifolds with contact type boundaries, Invent. Math. {\bf 103} (1991), 651-671.

\bibitem{medina-revoy-lattice85} \textsc{A. Medina and P. Revoy},  Les groupes oscillateurs et leurs r\`eseaux. Manuscripta Math. {\bf 52} (1985), 81-95.

\bibitem{medina-revoy-lattice07} \textsc{A. Medina and  P. Revoy}, Lattices in symplectic Lie groups. J. Lie Theory {\bf 17} (2007), 27-39.

\bibitem{milnor} \textsc{J. Milnor},  Curvatures of left invariant metrics on Lie groups,
Advances in Math. {\bf 21} (1976), 293--329.

\bibitem{mostow} \textsc{G. Mostow},  Factor spaces of solvable groups Ann. Math. {\bf 60} (1954), no. 1, 1-27.

\bibitem{ov1} \textsc{A. L. Onishchik} (Ed.),  Lie Groups and Lie Algebras I:  Foundations of Lie Theory and Lie Transformation Groups, Encyclopedia of Mathematical Sciences {\bf 20}, Springer-Verlag Berlin, Heidelberg, 1993.

\bibitem{ov2} \textsc{A. L. Onishchik and E. B. Vinberg} (Eds.),  Lie Groups and Lie Algebras II:  Discrete Subgroups of Lie Groups and Cohomologies of Lie Groups and Lie Algebras, Encyclopedia of Mathematical Sciences {\bf 21}, Springer-Verlag Berlin, Heidelberg, 2000.

\bibitem{raghu} \textsc{M. S. Raghunathan},  Discrete Subgroups of Lie Groups, Springer-Verlag, Berlin, Heidelberg, 1972.

\bibitem{sawai} \textsc{H. Sawai},  A construction of lattices on certain solvable Lie groups, Topology and Its Applications {\bf 154} (2007), 3125-3134.



\bibitem{sy}\textsc{H. Sawai and T. Yamada},  Lattices on Benson-Gordon type solvable Lie groups, Topology and Its Applications {\bf 149} (2005), 85-95.

\bibitem{wang} \textsc{H.-C. Wang},  Discrete subgroups of solvable Lie groups, Ann. Math. (1956), 1-19.

\end{thebibliography}
\end{document}